\documentclass[preprint,12pt]{elsarticle}

\usepackage{graphicx}
\usepackage{amssymb}
\usepackage{amsmath}
\biboptions{comma,square,sort&compress}

\journal{jcp}

\begin{document}

\begin{frontmatter}


%% \tnotetext[label1]{tnotetext}
%% \author{Mihail Semenov\corref{cor1}\fnref{label2}}
 \author{Mihail Semenov}
 \ead{sme@tpu.ru}
 \ead[url]{http://portal.tpu.ru/SHARED/s/SME/}
 %%\fntext[label2]{fntext}
 %%\cortext[cor1]{cortext}
 \address{30, Lenin Avenue, Tomsk, 634050, Russia, Tomsk Polytechnic University\fnref{label3}}
 %%\fntext[label3]{fntext}

\title{Analyzing the absolute stability region of implicit methods of solving ODEs}


%%\author{Mihail Semenov}
%%\address{Tomsk Polytechnic University}

\begin{abstract}
Review of implicit methods of integrating system of stiff ordinary
differential equations is presented. Defines and graphically
presents absolute stability region for Gears methods (backward
differentiation formula) used to solve system of stiff ordinary
differential equations. Recommendations for selecting the order of
Gears method are given.
\end{abstract}


\end{frontmatter}

% \linenumbers

\section{Introduction}\label{}

Ordinary differential equations (ODE) are widely used for
modelling real processes. Practice shows that the initial value
problem (Cauchy problem) for ODE systems can be attributed to the
following types: soft, stiff, ill-conditioned and rapidly
oscillating. Each type is connected with specific demands to
integrating methods. Stiff systems can be exemplified by problems
of chemical kinetics \cite{Kahaner89, Forsite80}, nonstationary
processes in complex electric circuits \cite{Hairer96, Butcher08},
systems emerging while solving equations of heat conduction and
diffusion \cite{Kalitkin95}, movement of celestial bodies and
satellites \cite{Gear87}%%, cellular reactions in the heart
%%\cite{MacLachlan07}
, plasticity physics \cite{Semenov07}, etc.

The numerical solution of ODE systems is accompanied by problems
due to the fact that at the modelling of a complex physical
process speeds of local processes may vary significantly, while
variables systems may be of different orders or change within the
interval of integration by orders of magnitude \cite{Semenov07,
Kolupaeva10}.

Besides, in physical experiments a boundary layer may be observed
which is characterized by quick changes in the object of research.
Such a layer does not necessarily emerge at the beginning of the
experiment, but only when some "controlling" parameter sharply
changes or achieves a certain critical value. For the numerical
solution of such problems the numerical method has to be chosen
very carefully \cite{Kahaner89, Forsite80, Cohen95, Enright87,
Dekker84, Rosenbrock68, Babenko79}.

The purpose of this paper is to report on a continuing research
effort aimed at the use and development of numerical methods in
computer program for  plasticity physics. The previous results
were reported in \cite{Kolupaeva10} while using the same
mathematical basis for solving stiff ODEs. The paper is organized
as follows: the review of definitions stiff ODEs presented in
Section \ref{sec2}. The approaches to finding a numerical solution
of the stiff ODEs are discussed in Section 3. In Section 4,
defines and graphically presents absolute stability region for
backward differentiation formula and recommendations for various
types of ODEs are given.

\section{Stiff ordinary differential equations}
\label{sec2} Interest for stiff systems appeared at the beginning
of the 20th century, initially in radio engineering (van der Pol
problem, 1920 \cite{Hairer96, Butcher08, Hull72, Shampine79,
Mazzia03, Brugnano09}). Then there was a new wave of interest in
the middle 1950s with application in studying equations of
chemical kinetics, movement of celestial bodies \cite{Gear87},
which contained both very slowly and very rapidly changing
components. Methods like Runge-Kutta, which had been considered as
highly reliable, produced mistakes in solving such problems.

One of the first attempts to give a definition of stiff systems
was made by C.F.~Curtiss and J.O.~Hirschfelder in 1952. They
proposed the following interpretation: stiff equations are
equations where certain implicit methods perform better, than
using classical explicit ones like Euler or Adams methods
\cite{Curtiss52}.

There are a number of interpretations of stiffness, each of which
reflects certain aspects of the numerical solution (e.g.
impossibility of using explicit methods of integration
\cite{Hairer87}, presence of rapidly damped disturbances
\cite{Dekker84, Ascher98}, large Lipschitz constants or
logarithmic norms of matrices \cite{Soderlind06, Higham93}, big
difference between eigenvalues of Jacobian matrix \cite{Kahaner89,
Forsite80, Ver02, Chua80}, fullness of Jacobian matrix
\cite{Nejad05}, a priori fixed sign of the solution
\cite{Bertolazzi96}, number of transient phases \cite{Novati03},
etc.) In some applications an important factor influencing the
behavior of the numerical solution is the order of the ODE system
\cite{Weiner97} or the presence of a boundary level \cite{Chua80,
Rakit79}, in others -- only the limit behavior (gradual change) at
wide intervals of integration. It is often unclear whether
stiffness is attributed to a particular solution or the problem in
general.

In papers \cite{Ortega86, Kahaner89, Forsite80, Hairer96,
Babenko79, Ver02, Chua80} authors find it difficult to clearly
define a stiff system ODE because of its complex character, so
they present a working description of a stiff problem. This is a
problem modelling a physical process, components of which possess
incommensurable characteristic times, or a process characteristic
time (reciprocal quantities of Jacobian eigenvalues) of which is
much smaller than the interval of integration.

In 1970s L.F.~Shampine and C.W.~Gear, who had gained a wide
experience of computing experiments with systems having components
of the decision vector of different orders, offered their own
definition of a stiff ODE system: the initial value problem for
ODEs is stiff if the Jacobian $J_{i,j}=\partial f_i/\partial y_j$,
$i$, $j=1, \ldots , N$ has at least one eigenvalue, for which real
part is negative with high modulus, while the solution within the
major part of the interval of integration changes slowly
\cite{Shampine79, Gear71a}.

In papers \cite{Ver02, Rakit79} ODE system is called stiff if real
components of all the eigenvalues of Jacobian are negative, i.e.
$Re(\lambda_i)<0$, $i=1, \ldots, N$ (the system is asymptotically
stable) and the ratio
$$
s = \frac
    {max\{|Re(\lambda_i)|, i = 1, 2, \ldots, N\}} %%\over
    {min\{|Re(\lambda_i)|, i = 1, 2, \ldots, N\}}
$$
is large. The parameter $s$ is called the stiffness ratio.

The problem of defining a stiff system is that for the stiffness
ratio the boundary value where it becomes big is not given. A
system of equations can be regarded as stiff if the stiffness
ratio $s$ exceeds 10, but in numerous applied problems this
parameter reaches $10^6$ and higher \cite{Petcu00}. In paper
\cite{Gridin03} a concept of superstiffness is introduced, when
the stiffness ratio reaches $10^6 \ldots 10^{12}$.

It must be pointed out that there are no simple methods for
evaluating stiffness, so numerical methods working without
stiffness testing are necessary.

\section{Approaches to finding a numerical solution of the stiff ODEs}

Let us regard a Cauchy problem for the ODEs of the first order,
which can be presented as follows:
\begin{equation}\label{eq1}
 Y'=F(x,Y), x\in[x_0, b], Y(x_0)=Y_0.
\end{equation}

For numerical integrating of system (\ref{eq1}) methods using a
linear combination of the decision vector $Y_1$, $Y_2$, \ldots ,
$Y_n$, \ldots and its derivatives in some sequence of independent
variable $x_1$, $x_2$, \ldots , $x_n$, \ldots are widely applied.
Such methods are called linear methods \cite{Kahaner89, Butcher08,
Gear71a}.

Since the discovery of the stiffness in the development of
numerical methods for integrating stiff systems, the following
trends have appeared: the investigation of stiffness and the
establishment of the theoretical apparatus stability analysis
methods \cite{Hairer87, Hairer96}, design and enhancement of the
methods taking into account the specifics of the tasks
\cite{Vshivkov07, Rahunanthan10, Filippov04, Hojjati04,
Skvortsov09, Hundsdorfer09, Savcenco09} and the prospects for
parallel \cite{Podhaisky02, Bahi09}. The most complete overview of
the current numerical methods for solving stiff ODE systems with
an extensive bibliography is presented in the papers
\cite{Hairer96, Butcher08, Curtiss52, Rakit79, Wang06, Hall79}.

\subsection{Taylor series}

One of the easy ways to construct the solution to the system
(\ref{eq1}) at point, $x_{n+1}$, if it is known at point, $x_n$,
is a method based on the expansion of solutions $Y(x_{n+1})$ in a
Taylor series in the neighborhood of point, $x_n$:
$$
Y(x_{n+1})=Y(x_n)+hF(x_n,Y_n,h),
$$
where
$$
F(x_n,Y_n,h) = Y'(x_n) + hY''(x_n)/2! + h^2Y'''(x_n) / 3 !+ \ldots
.
$$
If this series is truncated at $q$-th term and replace $Y(x_n)$
with the approximate value of $Y_n$, thus an approximate formula
is obtained:

\begin{eqnarray}\label{eq2}
Y_{n+1}=Y_n+h(F(x_n,Y_n)+h F'(x_n,Y_n)/2!+\nonumber\\
\qquad {}+ h^2F''(x_n,Y_n)/3!+\ldots+h^qF^{(q)}(x_n,Y_n)/(q+1)!).
\end{eqnarray}
If $q=1$, the computational scheme of explicit Euler method
\cite{Kahaner89, Forsite80, Babenko79, Ver01, Ver02} is obtained
$$
Y_{n+1}=Y_n+hF(x_n,Y_n).
$$

Application of the formula (\ref{eq2}) is limited to only those
tasks which can easily calculate the higher-order derivatives of
the function $F(x, Y)$ of the right side of the system
(\ref{eq1}). Though, it is usually not so.

\subsection{Runge-Kutta methods}

S.~Runge (1895), K.~Heun (1900) and M.~Kutta (1901) put forward an
approach based on constructing of the formula for $Y_{n+1}$
\cite{Kahaner89, Forsite80, Butcher08}:
$$
Y(x_{n+1})=Y(x_n)+h\Phi(x_n,Y_n,h),
$$
where $h$ -- integration step. The function $\Phi(\cdot)$ is close
to $F(\cdot)$, but does not contain the derivatives from the right
side of the equation. Thus, the series of explicit and implicit
methods requiring $s$-stage calculation of the right side function
at each integration step are obtained ($s$-stage methods).

The formulas of these methods are ideally applicable for practical
calculations: they allow to change the integration step $h$
easily. Perhaps the most famous is the formula of the 4-th order
of 4-stage Runge-Kutta \cite{Butcher08}.

One of the major problems associated with the use of Runge-Kutta
methods (in fact, almost of all the explicit methods) lies in
choosing the size of the integration step $h$, which provides the
stability of the computational scheme \cite{Kahaner89, Forsite80,
Babenko79, Ver01, Ver02, Chua80, Arushanyan90}. Nevertheless, even
nowadays, the explicit adaptive methods for solving stiff ODEs
\cite{Skvortsov09} are developed and widely used.

\subsection{Backward differentiation formulas}

These stiff tasks have made the implicit computational scheme
particularly attractive and led to the development of such
implicit methods which do not involve calculations based on the
size of the integration step \cite{Curtiss52, Chua80, Gear71a,
Gear71b, Gear71c, Nordsieck62, Miranker75, Miranker81}. The most
common among them are the methods of Adams-Moulton and "backward
differentiation formulas" (more commonly known as Gear method).
Having got the approximation to the solution at points, $x_1$,
$x_2$, \ldots, $x_n$, it is possible to find solutions at point,
$x_{n+1}$.

The computational schemes of the Adams-Moulton implicit methods
take the following form \cite{Kahaner89, Forsite80, Babenko79,
Ver01, Ver02, Chua80, Arushanyan90}:

\begin{equation}\label{eq3}
Y_{n+1}=Y_n+h\sum_{i=0}^q \beta_iF(x_{n-i+1}, Y_{n-i+1}),
\end{equation}
where $q$ determines the order of the method, the constants
$\beta_i$, $i=0, 1,\ldots, q$ correspond to the chosen order of
the method \cite{Butcher08, Chua80,  Gear71a}. Moreover, the
implicit Euler's method (the first order) and the trapezoidal
method (the second order) are the special cases of the last
computational scheme~(\ref{eq3}) where $q = 0$ and $q = 1$.

The construction of the multistep methods is based on the
polynomial of the degree $q$. The approximate value of the
solution $Y(x)$ at the point, $x_{n+1}$ appears as a linear
combination of the several approximate values of the solution and
its derivative at this and the previous $q$ points. Obviously, the
use of the multistep formulas requires the calculations of the $q$
units of the initial values $Y_1$, $Y_2$, \ldots, $Y_n$. The
accuracy of setting these $q$ units should not be less precise
than the accuracy of the formula. Thus, the polynomial can be
represented \cite{Kahaner89, Forsite80, Babenko79, Ver01, Ver02,
Chua80, Arushanyan90} as the following formula:

\begin{equation}\label{eq4}
Y_{n+1}=\sum_{i=1}^q \alpha_i Y_{n-i+1} + h\sum_{j=0}^q
\beta_jF(x_{n-j+1}, Y_{n-j+1}).
\end{equation}
Some constants $\alpha_i$ and $\beta_i$ in equation~(\ref{eq4})
can take zero values. When $\beta_1=\\=\beta_2=\ldots=\beta_q=0$,
it is possible to construct backward differentiation formulas.

\section{The region of implicit methods stability}

The first studies on the stability of multistep methods refer to
the researches of G.~Dahlquist \cite{Dahlquist56, Dahlquist63}.
According to the definitions above, the stiff ODEs require the
stability of numerical methods used to solve them. When getting
the asymptotically stable solution of the stiff Cauchy problem, an
error of the difference method should non-increase under any step,
i.e. the method should be absolutely stable. Current reviews of
the stability regions of multistep methods can be found in
\cite{Butcher08, Hairer96}.

The researchers \cite{Ver02, Ver01, Samarskiy89} give a more clear
definition to the stability method through the model first-order
equation
\begin{equation}\label{eq5}
y'= \lambda y, y(x_0)=y_0.
\end{equation}

The general solution of the equation~(\ref{eq5}) takes the form of
$$
y=C\cdot exp\{\lambda x\},
$$
where $C$ -- a constant, the solution of the %%corresponding Cauchy
%%problem with the initial conditions $y(x_0)=y_0$
equation~(\ref{eq5}) -- a function $y=y_0\cdot
exp\{\lambda(x-x_0)\}$ -- tends to zero, if $Re(\lambda)<0$ and it
infinitely grows in absolute value for $Re(\lambda)>0$, where
$\lambda$ is a complex number. Further, the concept of absolute
stability and the study of the absolute stability of the numerical
methods are considered for the model equation~(\ref{eq5}).

%%The research \cite{Hall79}[45] proves the validity of the model
%%equation~(\ref{eq5}).
Here, the stability domain of the multistep method~(\ref{eq4}) of
solving the initial value problem~(\ref{eq5}) is defined as the
set of points of complex numbers plane defined by the complex
variable $\sigma = h\lambda$. For $\sigma=h\lambda$, this method
applied to the model equation (\ref{eq5}) is stable, i.e. it
provides non-increase of an error \cite{Ver02, Ver01,
Samarskiy89}.

To determine the region of the implicit methods stability
(\ref{eq4}), the characteristic (complex) polynomial is used:
$$
\begin{array}{l}
P(z)=(z^q - \alpha_1z^{q-1}- \alpha_2z^{q-2}-\ldots-\alpha_q)+\\
+\sigma(\beta_0z^q + \beta_1z^{q-1}+
\beta_2z^{q-2}+\ldots+\beta_q)=0,
\end{array}
$$
it may be represented in the form of \cite{Butcher08, Hairer96,
Chua80, Gear71a, Gear87}:
\begin{equation}\label{eq6}
\sigma(\theta)=-\frac{\sum\limits_{k=0}^q\alpha_k e^{i
(q-k)\theta}} {\beta_0 e^{i q\theta}},
\end{equation}
%%\begin{equation}\label{eq6}
%%\sigma(\theta)=-\frac{\left(%
%%\begin{array}{c}
%%  - e^{i q\theta} + \alpha_1 e^{i (q-1)\theta}+\alpha_2 e^{i (q-2)\theta}+ \\
%%  +\ldots+\alpha_{q-1} e^{i \theta}+\alpha_q \\
%%\end{array}%
%%\right)}{\beta_0 e^{i q\theta}},
%%\end{equation}
where $\alpha_k$, $\beta_k$ -- coefficients
of the method ($k$=0, 1, 2,\ldots, $q$;
$\beta_1=\beta_2=\ldots=\\=\beta_q=0$), $q$ -- the order of the
method, $i$ -- a unit imaginary number, $z=e^{i\theta}$ -- an
imaginary number, $0\leq\theta\leq2\pi$. To determine the region
of the method absolute stability for the given value $\sigma=
\sigma_0$, it is necessary to find the solution of the
equation~(\ref{eq6}) relatively to $\theta$.

A set of points generated by the equation~(\ref{eq6}) corresponds
to a geometrical locus of points of single radicals
$\Gamma_\sigma$, for which $|z|=1$ is true. The region of absolute
stability of implicit methods~(\ref{eq4}) is considered to be an
external area $\Gamma_\sigma$ since at
$|\sigma|=h|\lambda|\rightarrow\infty$ implicit methods
(\ref{eq4}) are stable \cite{Butcher08, Hairer96, Chua80,
Gear71a}.

\begin{figure}
\center
\includegraphics[width=0.9\textwidth]{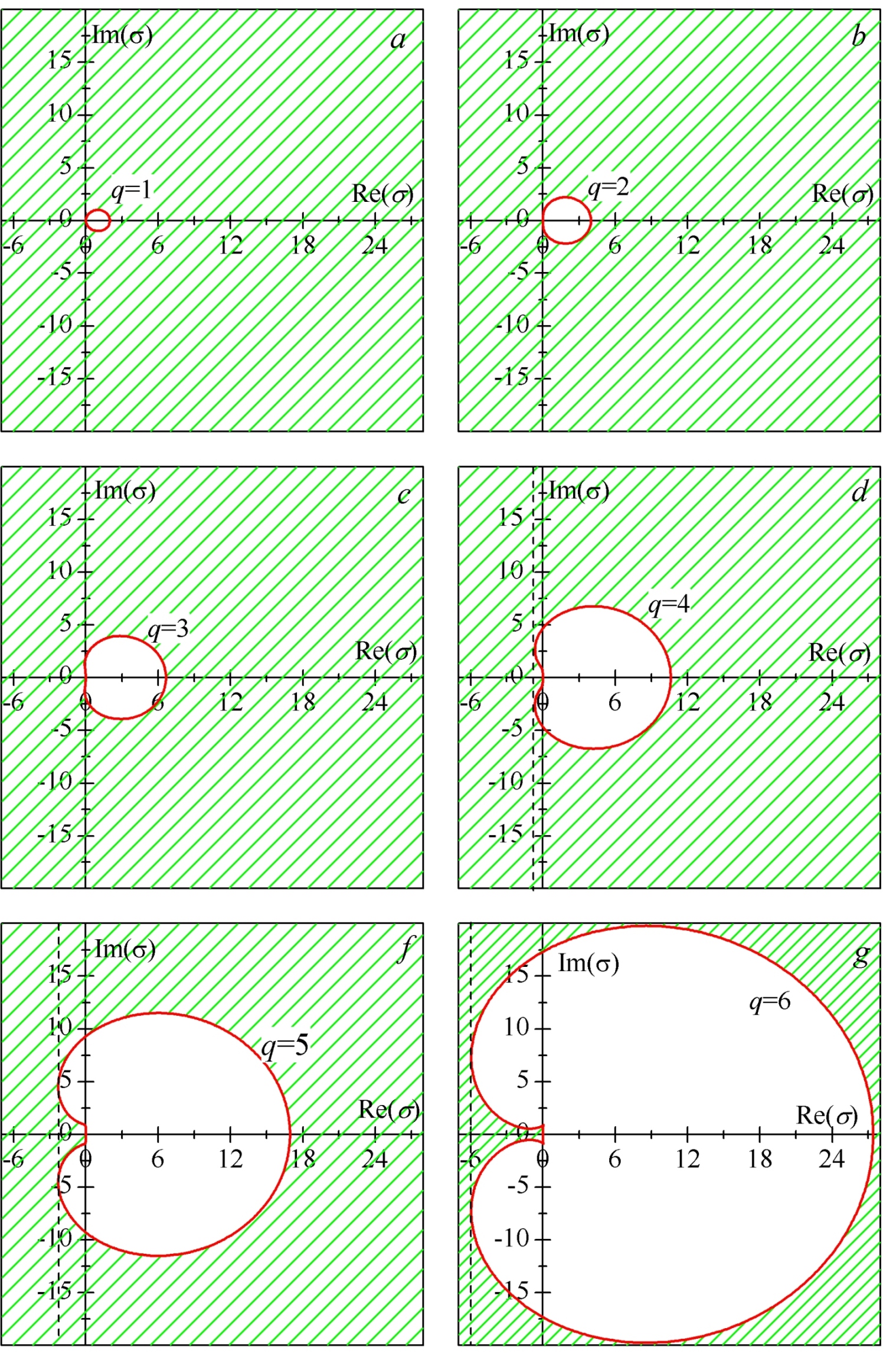}
\caption{Region of absolute stability for Gear methods of
$q$-order: a) first (Euler method); b) second; c) third; d) forth;
e) fifth; f) sixth.}\label{fig:1}
\end{figure}

To determine a geometrical locus of points described in the
equation~(6), a computer algebra system MathCAD is used. As a
result, a geometrical locus of single radicals $\Gamma_\sigma$ is
obtained (fig.~\ref{fig:1}). The region of absolute stability of
methods is considered to be an external region $\Gamma_\sigma$
(shaded).

The study region of stability ($\Gamma_\sigma$ is a simple closed
curve, fig.~\ref{fig:1}) shows that implicit methods~(\ref{eq4})
from 1 to 6 order inclusive are stiffly stable (first introduced
by \cite{Gear71a}). The stiffness property of the implicit
method~(\ref{eq4}) from 1 to 6 order inclusive is attained at
various values of a real number  $\delta\leq 0$ (fig.~\ref{fig:1},
$c-g$ marked with a dotted line). Particularly, for the methods of
the first and second order it is $\delta = - 0{.}$; for the third
order $\delta = - 0{.}1$; for the fourth order $\delta = - 0{.}7$;
for the fifth order $\delta =  -2{.}4$, for the sixth order
$\delta = -6{.}1$.

For Gear's method of the seventh order, the equation~(\ref{eq6})
will be:
$$%%\begin{equation}
\sigma(\theta)=-\frac{\left(%
\begin{array}{c}
  - e^{i 7\theta} + \frac{980}{363} e^{i 6\theta}-\frac{490}{121} e^{i 5\theta}-\frac{4900}{1089} e^{i 4\theta}- \\
  -\frac{1225}{363}e^{i 3\theta}+\frac{196}{121} e^{i2\theta}-\frac{490}{1089} e^{i\theta}+\frac{20}{363} \\
\end{array}%
\right)}{\frac{140}{369} e^{i 7\theta}}.
$$%%\end{equation}

\begin{figure}
\center
\includegraphics[width=0.90\textwidth]
{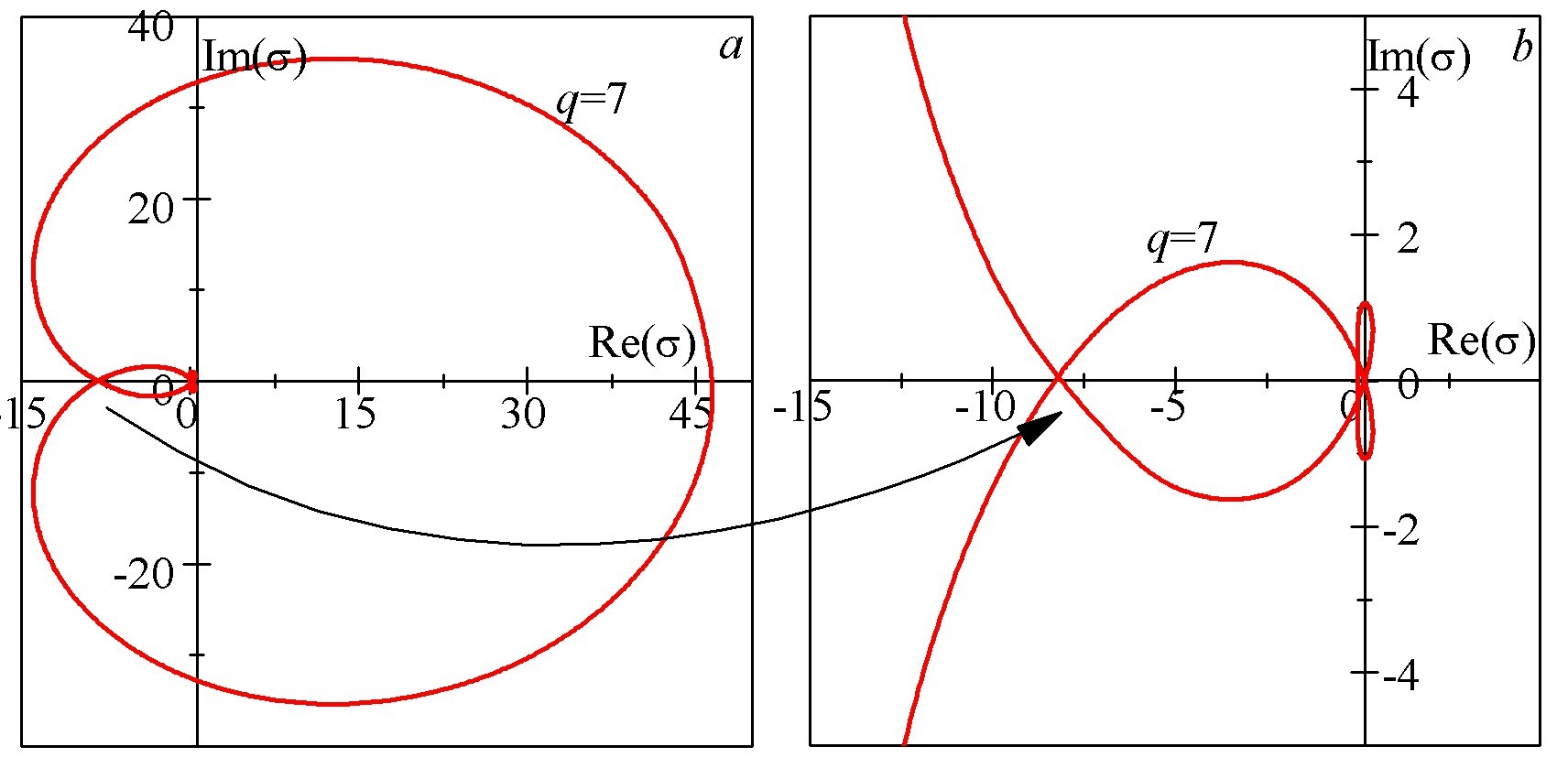}
\caption{The geometrical locus of points of single radicals
$\Gamma_\sigma$ for Gear method of the seventh order.}
\end{figure}\label{fig:2}

Finding a solution to $\sigma(\theta)$ relative to $\theta$, it
turns out that Gear method of the seventh order does not meet the
requirements of stiff stability (fig.~2)\cite{Chua80, Gear71a}. At
origin of coordinates and at $Re(\lambda)\approx-8$, there are
intersection points of $\Gamma_\sigma$.

The results obtained for the equation~(\ref{eq5}) can be
distributed on the ODEs \cite{Hall79}. In the case of an
autonomous system of ODE of $Y'=\textbf{A}Y$ type, where
\textbf{A} is a constant matrix, it becomes possible to transact
matching Jordan's form and proceed to look for a solution to the
system of ODE of $Z'=\textbf{J}Z$ type, where
$\textbf{J}=\textbf{T}^{-1}\textbf{AT}=diag\{\lambda_1, \lambda_2,
\ldots,\lambda_n\}$, $\lambda_i$ -- eigenvalues of matrix
\textbf{A}, $i=1, 2, \ldots, N$, $Y=\textbf{T}Z$,
$Z=\textbf{T}^{-1}Y$. The matrix \textbf{T} is composed of
eigenvectors of the matrix \textbf{A}. Thus the initial system of
ODE decomposes into $n$ scalar equations, for which the solution
can be found and the above approach to the region of stability
determination applied \cite{Butcher08, Chua80, Gear71a}.

If the coefficients of the system $\textbf{Y}'=\textbf{A}(x)Y$ are
not constant, the check of the eigenvalues \textbf{A} at each
value $x$ becomes laborious to calculate \cite{Gear71a}. It should
be noted that operating with nonlinear systems of
$Y'=\textbf{A}Y+\\+G(x,y)$ type, the stability of the solution can
be provided only at the origin of coordinates, moreover the
stability can be broken for eigenvalues located on an imaginary
axis \cite{Hairer87}.

The overview of the alternative ways of defining stability regions
for the implicit methods can be seen in \cite{Butcher08, Hairer87,
Hairer96}.

Implicit methods (Gear methods) can be applied for the calculation
of a big category of the stiff ODEs. In this case decrease of
stepsize (to the minimum possible) doesn't always let us adapt to
local solution and decrease volume of computation with required
precision. Optimal strategy of using multistep methods implies
availability of order autocontrol (from 1 to 6) and stepsize.

%%
%%It requires the calculation of local methodic error and effective
%%(in terms of computer storage saving) computation, storage and
%%reverse engineering of array variables preceding values of the
%%solution.

\section{Conclusion}

The area of the numerical methods for solving of ODEs is one of
the most well-investigated topics in the mathematical literature.
A number of techniques and solvers have been suggested,
development, and described, but a clear definition of stiffness
has not been provided, so the working definition of stiffness is
still topical.

As a rule in most solvers for ordinary differential equations the
explicit first order method (Euler method) or a second order
method (trapezoidal method) is applied. Implicit Gear methods
(backward differentiation formulas) are stiff from 1 to 6 order
inclusive, so for the acceleration of the integration process of
ordinary differential equations increasing order could be applied.

The results of the calculations let us define absolute stability
regions for the implicit methods where changing of integration
step over wide region when computational stability of the method
is constant.

\section{Acknowledgments}

The author is grateful to Dr. S.N.~Kolupaeva at Tomsk State
University of Architecture and Building, Russia and
Dr.~J.C.~Butcher at The University of Auckland, New Zealand who
have been influence in the development of the ideas presented in
this paper.



%% References


\bibliographystyle{elsarticle-num}
\bibliography{mybib}






% \end{thebibliography}

\end{document}